\documentclass[10pt,a4paper]{article}
\usepackage[utf8]{inputenc}
\usepackage[T1]{fontenc}
\usepackage{amsmath}
\usepackage{amssymb}
\usepackage{graphicx}

\newcommand{\ZZ}{{\mathbb Z}}

\newtheorem{cor}{\noindent\rm\bf Corollary}[section]
\newtheorem{thm}{\noindent\rm\bf Theorem}[section]

\newtheorem{defn}{\noindent\rm\bf Definition}[section]

\newtheorem{lem}{\noindent\rm\bf Lemma}[section]

\begin{document}
	
\title{The unimodular equivalence of sublattices in an $n$-dimensional lattice}
\author{Shikui Shang}
\date{}
\maketitle


\vspace{2mm}

\abstract{	
In this paper, we study the unimodular equivalence of sublattices in an $n$-dimensional lattice. A recursive procedure is given to compute the cardinalities of the unimodular equivalent classes with the indices which are powers of a prime $p$. We also show that these are integral polynomials in $p$. When $n=2$, the explicit formulae of the cardinalities are presented depending on the prime decomposition of the index
$m$. We also give an explicit formula on the number of co-cyclic sublattices with a fixed index $m$, which consist a unimodular equivalent class of sublattices.
}

{\bf 2020 Mathematics Subject Classification:} primary 11H06; secondary 20E07, 52C05, 52C07. 

{\bf Keywords and phrases:} Lattice; Sublattice; the unimodular equivalence; co-cyclic sublattices.

\section{Introduction}	

Let $\Lambda$ be an $n$-dimensional lattice in an Euclidean space. A subset $\Lambda^*$ of $\Lambda$ is called a sublattice of $\Lambda$ if it is also an $n$-dimensional lattice.
Define the positive integer $m=[\Lambda:\Lambda^*]$ as the index of $\Lambda^*$ in $\Lambda$ and therefore
$$\det(\Lambda^*)=m\cdot\det{\Lambda}.$$

In \cite{Z}, the concept of unimodular equivalence of sublattices is introduced. Two sublattices $\Lambda^*$ and $\Lambda^\bullet$ of $\Lambda$ are called unimodular equivalent if there exists a linear transformation $\sigma$ such that
$$\Lambda=\sigma(\Lambda)\ \ \text{  and   }\ \ \Lambda^\bullet=\sigma(\Lambda^*).$$ 
Obviously, the index is an invariant of unimodular equivalent sublattices.

The reason that the equivalence is called unimodular is that the linear transformation $\sigma$ corresponds to a unimodular matrix in $GL_n(\ZZ)$.

The structures and representations of sublattices have been studied by Minkowski, Siegel, Cassels, Hlawka, Rogers, Schmidt and Gruber.  Particular, sublattices have been studied in \cite{B}-\cite{F},\cite{Sc1},\cite{Sc2}.

Fix an $n$-dimensional lattice $\Lambda$ and an integer $m>0$.
Let $f_n(m)$ denote the number
of the different sublattices of $\Lambda$ with index $m$ and let $g_n(m)$ denote the number of the different
sublattice classes of $\Lambda$ with index m under unimodular equivalence.

In 1945--1946, C.L.Siegel gave the first upper bound for $f_n(m)$, namely,
$$f_n(m)\leq m^{n^2},$$
in his celebrated lecture notes \cite{Si}. 

Although counting the subgroups of a group is a classic topic in algebra (see the classic books \cite{LS}, \cite{Se} and papers such as \cite{GLP}, \cite{GSS}), explicit formulae for $f_n(m)$ were achieved only in 1997(see \cite{G}).  In 1997, Baake \cite{B} also deduced the following formula of $f_n(m)$ based on a
recursion,
\begin{equation}f_n(m)=\sum_{d_1\cdots d_n=m}d_1^0d_2^1\cdots d_n^{n-1}\tag{1.1}.\end{equation} 

 Recently, C.Zong gave a systematic treatment on  $f_n(m)$ and $g_n(m)$ in \cite{Z}. He proved that if the integer $m$ has the prime decomposition $m=p_1^{r_1}\cdots p_l^{r_l}$, then
\begin{equation}f_n(m)=\prod_{i=1}^l\prod_{j=1}^{r_i}\frac{p_i^{j+n-1}-1}{p_i^j-1}=\prod_{i=1}^l\prod_{j=1}^{n-1}\frac{p_i^{j+r_i}-1}{p_i^j-1}\tag{1.2}\end{equation}
(also see \cite{G}) and
\begin{equation}g_n(m)=\prod_{i=1}^lp_n(r_i),\tag{1.3}\end{equation}
where $p_n(k)$ denote the number of partitions of $k$ into  $n$ non-negative
parts. More in detail, define the set of partitions of $k$ into  $n$ non-negative parts as
$${\mathcal P}(n,k)=\{(\alpha_1,\cdots,\alpha_n)\in\ZZ^n\ |\ 0\leq\alpha_1\leq\cdots\leq\alpha_n, \sum_{i=1}^n\alpha_i=k\}.$$
then $p_n(k)$ is the cardinality of ${\mathcal P}(n,k)$.

\

It is natural to ask how the size of a given unimodular equivalent class is and how sublattices are distributed in them. In this paper, we study these problems and give some results for them.

First, we observe that $f_n(m)$ and $g_n(m)$ are multiplicative on $m$, which make us focus on the cases where $m=p^r$ is a power of prime. 
Next, $f_n(m)$(resp. $g_n(m)$) is the number of integral matrices of Smith(resp. Hermite) normal form with determinants $m$. To determine the unimodular equivalent class of a given sublattice is the same thing that to compute the Smith normal form of an integral matrix of Hermite normal form. We can prove that the cardinalities of unimodular equivalent classes are also multiplicative. These basics are shown detailedly in Section 2.

We give some results (Theorem \ref{thm:3.1} and \ref{thm:3.2}) for small $n$($=2,3$) in Section 3. In Section 4, a recursive formula of the number of sublattices in a unimodular equivalent class is given when the indices of sublattices are powers of a prime $p$ (Theorem \ref{thm:4.2}). Using this formula, we show that all these numbers are polynomials in $p$ (Corollary \ref{cor:4.1}). 

A sublattice $\Lambda^*$ of $\Lambda$ is called co-cyclic if the quotient $\Lambda/\Lambda^*$ is a finite cyclic group. In other words, $\Lambda=\Lambda^*+\ZZ b$ for a lattice vector $b\in\Lambda$. In \cite{NS}, P. Nguyen and I. Shparlinski first count the number of co-cyclic sublattices and  gave an asymptotic formula of the number $N_n(V)$ of co-cyclic sublattices with index at most $V$. That is,
\begin{equation}N_n(V)=\frac{\vartheta_n}{n}V^n+O(V^{n-1+o(1)}),\tag{1.4}\end{equation}
where $\vartheta_n$ is an Euler product given by $$\vartheta_n=\prod_p\left(1+\frac{p^{n-1}-1}{p^{n+1}-p^n}\right).$$
They also gave an asymptotic formula on the number of co-cyclic sublattices with square-free index. Such co-cyclic sublattices have been previously studied from a complexity point of view in \cite{PS} and \cite{T}.

In the last section of our paper, we prove that all co-cyclic sublattices with the same index $m$ are unimodular equivalent(Theorem \ref{thm:5.1}). The explicit formulae on the number of co-cyclic sublattices are given successively when $m$ are power-of-primes, square-free and in general(See Theorem \ref{thm:5.1}, Corollary \ref{cor:5.1} and  Corollary \ref{cor:5.2}). Finally, We show another aspect on the fact that most sublattices are co-cyclic mentioned in \cite{NS}. 

\
	
\section{Basics}	
\subsection{The multiplicativity of $f_n(m)$ and $g_n(m)$}
	
In $(1.2)$ and $(1.3)$, we have seen that both $f_n(m)$ and $g_n(m)$ are multiplicative, i.e., $$f_n(1)=g_n(1)=1$$ and $$f_n(m_1m_2)=f_n(m_1)f_n(m_2), g_n(m_1m_2)=g_n(m_1)g_n(m_2)$$ wherever the integers $m_1$ and $m_2$ are coprime. This result can be obtained by the following theorem.
	
\begin {thm}\label{thm:2.1} Let $\Lambda$ be an $n$-dimensional lattice and $\Lambda^*$ a sublattice of $\Lambda$ with index $m_1m_2$ and $\gcd(m_1,m_2)=1$. Then, there exists a unique sublattice $\Lambda'$ of $\Lambda$ which contains $\Lambda^*$ as its sublattice satisfying $[\Lambda:\Lambda']=m_1$ and $[\Lambda':\Lambda^*]=m_2$.
\end{thm}	
{\bf Proof.} The sublattices of $\Lambda$ containing $\Lambda^*$ correspond to the subgroups of the quotient $\Lambda/\Lambda^*$, which is a finite abelian group. By Sylow $p$-subgroup theorem or the structure theorem of finitely generated abelian groups, we obtain the result.  $\Box$

{\bf Remark. }When $n=2$, the above theorem is also the reason of the multiplicativity of the Hecke operators $T_m$ (see, for example, Ch. 7 of \cite{Se}).
	
Hence, if we konw the prime factorization $m=p_1^{r_1}\cdots p_l^{r_l}$ of $m$, then
$$f_n(m)=f_n(p_1^{r_1})\cdots f_n(p_l^{r_l}),\ \ g_n(m)=g_n(p_1^{r_1})\cdots g_n(p_l^{r_l}),$$
which tell us that the values of $f_n(m)$(resp. $g_n(m)$) can be determined by their values on the powers of primes.
	
\subsection{Hermite Normal form and Smith Normal form}

Fix a basis of a lattice $\Lambda$, there exists a one-to-one correspondence between the sublattices of $\Lambda$ with index $m$ and the integral matrices $H$ with determinants $m$ of Hermite normal form(HNF),
i.e., $$H=\begin{pmatrix}
		h_1 &  &   &   &   \\
		& h_2 &  & h_{ij} &   \\
		&    &  h_3 &  &   \\
		&   0&   & \ddots &  \\
		&   &   &   & h_n \\
	\end{pmatrix},\text{ with }h_1h_2\cdots h_n=m,
	$$
where $H$ is upper-triangular and $0\leq h_{ij}<h_j$ for all $1\leq i\leq j-1$.
	
Denote	$$\text{HNF}_n(m)=\{H\in{\mathbb Z}^{n\times n}\ |\ H \text{ is of HNF with }\det(H)=m \}.$$
We have that $f_n(m)$ is the cardinality of $\text{HNF}_n(m)$. 
	
Correspondingly, there exists a one-to-one correspondence between the unimodular
equivalent classes of sublattices with index $m$ and the integral matrices $S$ with $\det S=m$ of Smith normal form(SNF), i.e., $$S=\begin{pmatrix}
		d_1 &  &   &      \\
		& d_2 &  &     \\
		&   &    \ddots & \\
		&   &     & d_n \\
\end{pmatrix}, \text{with }d_1d_2\cdots d_n=m$$ where all $d_i$ are positive integers satisfying $d_1\mid d_2\mid\cdots\mid d_n$.

Set
$$\text{SNF}_n(m)=\{S\in{\mathbb Z}^{n\times n}\ |\ S \text{ is of SNF with }\det S=m \}.$$
Then, $g_n(m)$ is the cardinality of $\text{SNF}_n(m)$. 

\begin{thm}\label{thm:2.2} If $m$ has the prime factorization $m=p_1^{r_1}\cdots p_l^{r_l}$, there exist bijections
$$\Phi:\text{HNF}_n(m)\simeq \text{HNF}_n(p_1^{r_1})\times\cdots\times \text{HNF}_n(p_l^{r_l})$$
and 
$$\Psi:\text{SNF}_n(m)\simeq \text{SNF}_n(p_1^{r_1})\times\cdots\times \text{SNF}_n(p_l^{r_l}).$$\end{thm}
{\bf Proof. }Use Chinese Remainder Theorem. $\Box$

{\bf Remark. }This result gives another explanation for the theorem in above subsection.	

Notice that for a prime $p$ and an integer $r>0$, any diagonal element of a matrix in $\text{HNF}_n(p^r)$ or $\text{SNF}_n(p^r)$ is also a powers of $p$ with exponent less than or equal to $r$. Under the mapping $\Phi$, $H\mapsto(H_1,\cdots,H_l)$ and $H_k\in \text{HNF}_n(p_k^{r_k})$ is of HNF with diagonal elements of $p_k$'s powers. In other words, if $\Lambda^*$ is a sublattice
with coefficient matrix $H\in \text{HNF}_n(m)$ in a given basis of $\Lambda$, then $\Lambda^*$ is the intersection of sublattices $\Lambda_k^*$ of $\Lambda$ with coefficient matrices $H_k$ (w.r.t. the same basis) with diagonals $p_k^{\alpha_{ki}}$ for $(\alpha_{k1},\cdots,\alpha_{kn})\in{\mathcal P}(n,r_k), 1\leq k\leq l$. And, $\Lambda^*_k$ can be seen as the $p_k$-component of $\Lambda^*$ in this sense.

The similar discussion is also adaptive for the unimodular equivalent classes.

\subsection{The invariant factors of integral matrices}	

In this subsection, we give information of the SNFs for integer matrices.

Let $M$ be an $(m\times n)$-matrix with integer coefficients. For $1\leq k\leq \min\{m,n\}$, denote  the set of all submatrices of $M$ with rank $k$ by ${\mathcal M}_k$ and
$$D_k(M)=\gcd\{\det(M_k) \ |\ M_k\in{\mathcal M}_k\}.$$
And, let $d_1(M)=D_1(M), d_k(M)=\frac{D_k(M)}{D_{k-1}(M)}, 2\leq k\leq n$, where $d_i(M)$ are called the invariant factors of $M$ and the powers of primes in the factorization of $d_i(M)$ are called the elementary factors of $M$.
	
By the knowledge of linear algebra, 
\begin{thm}\label{thm:2.3} If $M$ is an integral square martrix with rank $n$,  then $$\begin{pmatrix}
		d_1(M) &  &   &      \\
		& d_2(M) &  &     \\
		&   &    \ddots & \\
		&   &     & d_n(M) \\
	\end{pmatrix}$$ is the SNF of $M$,
	i.e., there are two unimodular matrices $U,V$ such that $$M=U\begin{pmatrix}
		d_1(M) &  &   &      \\
		& d_2(M) &  &     \\
		&   &    \ddots & \\
		&   &     & d_n(M) \\
	\end{pmatrix}V. \ \ \Box$$\end{thm}
We see that the invariant factors and the elementary factors of the coefficient matrices of a sublattice are independent of the choice of lattices basis. In fact, if $d_i$ are invariant factors with respect to the sublattice $\Lambda^*$, then $\Lambda/\Lambda^*$ is isomorphic to $(\ZZ/d_1\ZZ)\times\cdots\times(\ZZ/d_n\ZZ)$ as abstract finite abelian groups.	Furthermore, two sublattices are unimodular equivalent if and only if their invariant factors coincide. Hence, we can talk the invariant factors (and the elementary factors) 
of a unimodular equivalent class of sublattices.

\begin{defn}\label{def:2.1} Let $\Lambda$ be an $n$-dimensional lattice. For positive integers $d_1\mid\cdots\mid d_n$, define $\mathcal S(d_1,\cdots,d_n)$ is the unimodular equivalent class of sublattices in $\Lambda$ with invariant factors $d_1,\cdots,d_n$. Moreover, denote the cardinality of ${\mathcal S}(d_1,\cdots,d_n)$ by $f(d_1,\cdots,d_n)$.\end{defn}

It is obvious that
$$f_n(m)=\sum_{d_1\mid\cdots\mid d_n, d_1d_2\cdots d_n=m}f(d_1,\cdots,d_n).$$

Using Theorem \ref{thm:2.2} and \ref{thm:2.3}, we have
\begin{cor}\label{cor:2.1}  There is a bijection between two sets 
$$\mathcal S(d_1,\cdots,d_n)\simeq\{H\in \text{HNF}_n(m)\ |\ d_k(H)=d_k \text{ for all }1\leq k\leq n\}.$$
Then,  $f(d_1,\cdots,d_n)$ is equal to the number of $H\in \text{HNF}_n(m)$ with the specified invariant factors $d_k(H)=d_k, 1\leq k\leq n$.  $\Box$\end{cor}

For a fixed prime $p$ and an index $\alpha=(\alpha_1,\cdots,\alpha_n)\in{\mathcal P}(n,k)$, denote the notions 
$${\mathcal S}_\alpha(p)={\mathcal S}(p^{\alpha_1},\cdots,p^{\alpha_n}),\ \ f_\alpha(p)=f(p^{\alpha_1},\cdots,p^{\alpha_n}).$$

We also have multiplicativity of prime components,
\begin{cor}\label{cor:2.2}   Assume $m=p_1^{r_1}\cdots p_l^{r_l}$ and $d_i=p_1^{\alpha_{1i}}\cdots p_l^{\alpha_{li}}$ for $1\leq i\leq n$ with $\alpha_k=(\alpha_{k1},\cdots,\alpha_{kn})\in{\mathcal P}(n,r_k)$ for all $1\leq k\leq l$. Then,
	$$\mathcal S(d_1,\cdots,d_n)\simeq\mathcal S_{\alpha_1}(p_1)\times\cdots\times{\mathcal S}_{\alpha_l}(p_l),$$
	where $p_k^{\alpha_{ki}}$ are the elementary factors of sublattices in $\mathcal S(d_1,\cdots,d_n)$. Moreover,
	$$f(d_1,\cdots,d_n)=f_{\alpha_1}(p_1)\cdots f_{\alpha_l}(p_l). \ \ \Box$$
\end{cor}

{\bf Remark. }(1) Notice that it is simple to deal with the gcd of powers of a single prime $p$. It is just the power of $p$ with minimal exponent. Then, some counting problems reducing to powers of primes may become handled. 

(2) Although in the algorithmic view, to determinate the SNF of a concrete integral matrix is not too difficult(use elementary transformations of matrices over $\ZZ$), our problem is whether there exsits a closed formular for $f(d_1,\cdots,d_n)$ as the formulae $(1.1)$ or $(1.2)$ for $f_n(m)$.

\

\section{The cases for $n=2,3$}

For the lattices with small dimension $n=2,3$, some explicit results of unimodular equivalence classes can be obtained by detailed analysis.

As above mentioned, we consider locally a fixed prime. For a prime $p$ and a nonzero integer $m$, set $\text{ord}_p(m)$ is the maximal positive integer $t$ such that $p^t|m$ and $\text{ord}_p(0)=\infty$.

For $n=2$, assume that $m=p^r, r>0$. Let $H\in\text{HNF}_2(p^r)$ be of form
$\begin{pmatrix}
	p^{r_1} & h_{12} \\
	0 & p^{r_2} \\
\end{pmatrix}
$ with $r_1+r_2=r$ and $0\leq h_{12}<p^{r_2}$. Set $t=\min\{r_1,r_2,\text{ord}_p(h_{12})\}$.
Then, the SNF of $H$ is
$\begin{pmatrix}
	p^t & 0 \\
	0 & p^{r-t} \\
\end{pmatrix}$. Notice that $r_2\leq \text{ord}_p(h_{12})$ if and only if $h_{12}=0$. 

Consider the unimodular equivalent class $\mathcal S_{(t,r-t)}(p)$ of sublattices with the invariant factors $p^t, p^{r-t}$ satisfying $2t\leq r$. We have

\begin{thm}\label{thm:3.1} For $m=p^r$, 
$$f_{(t,r-t)}(p)=\begin{cases}p^{r-2t}+p^{r-2t-1},&0\leq t<\frac{r}{2}\\1, &t=\frac{r}{2}\end{cases}.$$
In particular, when $r$ is even, the unimodular equivalent class ${\mathcal S}_{\left(\frac{r}{2},\frac{r}{2}\right)}(p)$ has only one sublattice.
\end{thm}
{\bf Proof. }By Theorem \ref{thm:2.1}, for $H=\begin{pmatrix}
	p^{r_1} & h_{12} \\
	0 & p^{r_2} \\
\end{pmatrix}
\in \text{HNF}_2(p^t)$ with SNF $\begin{pmatrix}
	p^t & 0 \\
	0 & p^{r-t} \\
\end{pmatrix}$, one of the following three disjoint conditions is satisfied

(i)  $t=r_2, r_2\leq r_1, r_2\leq\text{ord}_p(h_{12})$,

(ii) $t=\text{ord}_p(h_{12}),\text{ord}_p(h_{12})<r_2, \text{ord}_p(h_{12})\leq r_1$,

(iii) $t=r_1, r_1<r_2, r_1<\text{ord}_p(h_{12})$.

By directly counting, we have
\begin{align*}
	&f_{(t,r-t)}(p)=\sum_{t=r_2\leq r_1}1+\sum_{t\leq r_1,t<r_2}(p^{r_2-t}-p^{r_2-t-1})+\sum_{t=r_1<r_2}p^{r_2-r_1-1}\\
	=&1+\sum_{t<r_2\leq r-t}(p^{r_2-t}-p^{r_2-t-1})+(1-\delta_{t,\frac{r}{2}})p^{r-2t-1}\\
	=&p^{r-2t}+(1-\delta_{t,\frac{r}{2}})p^{r-2t-1}=\begin{cases}p^{r-2t}+p^{r-2t-1},&0\leq t<\frac{r}{2}\\1, &r\text{ is even and }t=\frac{r}{2}\end{cases},
\end{align*}
where $\delta_{m,n}$ is the standard Kronecker symbol.  $\Box$

Recall that 
\begin{align*}&f_2(p^r)=p^r+p^{r-1}+\cdots+1\\
=&\underbrace{(p^r+p^{r-1})}_{f_{(0,r)}(p)}+\underbrace{(p^{r-2}+p^{r-3})}_{f_{(1,r-1)}(p)}+\cdots.\end{align*}
Hence, when $n=2$, it is clear that how sublattices constitute unimodular equivalent classes by the proof of Theorem \ref{thm:3.1}.

Using the multiplicativity in Corollary \ref{cor:2.2}, we also have 
\begin{cor}\label{cor:3.1}  Asumme that $m=p_1^{r_1}\cdots p_l^{r_l}$, $m_1=p_1^{t_1}\cdots p_l^{t_l}$ and $m_2=p_1^{r_1-t_1}\cdots p_l^{r_l-t_l}$ with $2t_k\leq r_k, 1\leq k\leq l$. Then,
$$f(m_1,m_2)=\prod_{1\leq k\leq l, t_k\neq\frac{r_k}{2}}(p^{r_k-2t_k}+p^{r_k-2t_k-1}). \ \ \Box$$\end{cor}

\

For $n=3$,  we only give a result from HNF to its SNF as follows.

Assume that $H\in HNF_3(p^r)$ is of form
$\begin{pmatrix}
	p^{r_1} & h_{12} & h_{13} \\
	0 & p^{r_2} & h_{23} \\
	0 & 0 & p^{r_3} \\
\end{pmatrix}
$ with $r_1+r_2+r_3=r, 0\leq h_{12}<p^{r_2}$ and $0\leq h_{13}, h_{23}<p^{r_3}$. Set $u=\text{ord}_p(h_{12}h_{23}-p^{r_2}h_{13})$. Then,

\begin{thm}\label{thm:3.2} The sublattice of $3$-dimensional lattice with coefficient matrix $H$ belongs to $\mathcal S(p^s,p^{t-s},p^{r-t})$ for
\begin{align*}
	&s=\min\{r_1,r_2,r_3,\text{ord}_p(h_{12}),\text{ord}_p(h_{13}),\text{ord}_p(h_{23})\},\\
	&t=\min\{r_1+r_2,r_2+r_3,r_1+r_3,r_1+\text{ord}_p(h_{23}),r_3+\text{ord}_p(h_{12}),u\}.
\end{align*}\end{thm}
{\bf Proof. } By Corollary \ref{cor:2.1}, we need to determinate the invariant factors of $H$. In fact, $p^s=d_1(H)=D_1(H)$ is the gcd of all coefficients of $H$. Moreover,
there are $6$ submatrices with rank $2$ $$H_{11},H_{21},H_{22},H_{31},H_{32},H_{33}$$ of $H$ with nonzero determinants, where $H_{ij}$ is the rank $2$ submatrix of $H$ omitting the $i$-th row and $j$-th column.  We obtain the results. $\Box$

It seems that when $n>3$, it becomes too complicated to give general explicit formulae of invariant factors even for a given HNF. 

\
	
\section{A recursive formula of $f_\alpha(p)$}

In Section 3, we have seen that to approach closed formulae of $f_\alpha(p)$ directly is difficult for an arbitrary partition $\alpha\in{\mathcal P}(n,k)$  when $n$ is large. One question is whether there is a recursive formula to deduce $f_\alpha(p)$. For this purpose, we try to compute $f_\alpha(p)$ by induction on $n$ in this section.
	
In geometric view, consider a sublattice $\Lambda^*$ in $\Lambda$. We first choose a basis $\{b_1,\cdots,b_n\}$ of $\Lambda$ and set the $(n-1)$-dimensional lattice $\Lambda_{n-1}=\ZZ b_2+\cdots+\ZZ b_n$ with a sublattice
$\Lambda^*_{n-1}=\Lambda^*\cap\Lambda_{n-1}$.

There exist a positive integer $h_1$ and a lattice vector $b_1^*\in\Lambda^*$ such that $$\Lambda^*=\Lambda_{n-1}^*+\ZZ b_1^*, \text{\ and\ \ } b_1^\star=b_1^*-h_1b_1\in\Lambda_{n-1}.$$

The following lemma gives us a chance to descent the dimension for counting sublattices, which is easy to be proved.
\begin{lem}\label{lem:4.2}Two sublattices $\Lambda^*_1$ and $\Lambda^*_2$ in $\Lambda$ coincide with each other if and only if 
$$\Lambda^*_{1,n-1}=\Lambda^*_{2,n-1}, h_1(\Lambda^*_1)=h_1(\Lambda^*_2), b_1^\star(\Lambda^*_1)\equiv b_1^\star(\Lambda^*_2)(\bmod\ \Lambda^*_{1,n-1}).\ \Box$$
\end{lem}
	
Hence, to count sublattices in a prescribed subset is equivalent to count the	appropriate triples $(h_1,\Lambda^*_{n-1},b_1^\star)$ for a positive integer $h_1$, a sublattice $\Lambda^*_{n-1}$ in $\Lambda_{n-1}$ and $b_1^\star\in\Lambda_{n-1}/\Lambda^*_{n-1}$. 

\

Assume that an HNF is of form $$H=\begin{pmatrix}
	h_1 &  &   &   &   \\
	0  & h_2 &  & h_{ij} &   \\
	0 & 0   &  h_3 &  &   \\
	\vdots&   &   & \ddots &  \\
	0 &0   & \cdots  &   & h_n \\
\end{pmatrix}=\begin{pmatrix}
	h_1 & h_{12}\cdots h_{1n} \\
	0 & H' \\
\end{pmatrix},
$$
where $H'$ is the bottom-right-cornered $(n-1)\times(n-1)$ submatrix of $H$. And, set $H''=\begin{pmatrix}
	h_{12}\cdots h_{1n} \\
	H' \\
\end{pmatrix}\in\ZZ^{n\times(n-1)}$.

\begin{lem}\label{lem:4.1}For $1\leq k\leq n-1$, the invariant factors of $H$ are
	
	$$D_k(H)=\gcd(h_1D_{k-1}(H'),D_k(H'')).$$\end{lem}
{\bf Proof. }Take a rank $k$ submatrix $H_k$ of $H$. If $H_k$ contains the first row and column, we have $H_k=\begin{pmatrix}
	h_1 & * \\
	0 & H'_{k-1} \\
\end{pmatrix}$ for a rank $(k-1)$ submatrix $H'_{k-1}$ of $H'$.
If $H_k$ contains the first column but doesn't contain the first row, $H_k$ has a zero column and $\det H_k=0$. Finally, if $H_k$ doesn't contain the first column, then $H_k$ is a rank $k$ submatrix of $H''$. We obtain the result. $\Box$	

For a lattice vector $b\in\Lambda$, denote
$$\Lambda^*(b)=\Lambda^*+\ZZ b,$$
which is a sublattice of $\Lambda$ containing $\Lambda^*$. 

Lemma \ref{lem:4.1} tells us that the invariant factors of $\Lambda^*$ can be determined by $h_1$ and the invariant factors of two $(n-1)$-dimensional lattices, $\Lambda^*_{n-1}$ and $\Lambda^*_{n-1}(b^\star)$.
	
\begin{thm}\label{thm:4.1} For a prime $p$, let $\Lambda^*$ be the sublattice of $\Lambda$ with coefficient matrix $$H=\begin{pmatrix}
		p^{r_1} &  &   &   &   \\
		0  & p^{r_2} &  & h_{ij} &   \\
		0 & 0   &  p^{r_3} &  &   \\
		\vdots&   & \ddots  & \ddots &  \\
		0 &0   & \cdots  & 0  & p^{r_n} \\
	\end{pmatrix}$$ in HNF. Set $$b_1^\star=b^*_1-p^{r_1}b_1=\sum_{i=2}^nh_{1i}b_i.$$ Then, for a partition $\alpha\in{\mathcal P}(n,\sum_{i=1}^nr_i)$, $\Lambda^*\in\mathcal S_\alpha(p)$ if and only if $\Lambda_{n-1}^*\in\mathcal S_\beta(p)$ and $\Lambda_{n-1}^*(b_1^\star)\in\mathcal S_\gamma(p)$ where the two partitions $\beta=(\beta_1,\cdots,\beta_{n-1}),\gamma=(\gamma_1,\cdots\gamma_{n-1})$ satisfy the following equations
\begin{align}
&\sum_{i=1}^k\alpha_i=\min\{r_1+\sum_{i=1}^{k-1}\beta_i,\sum_{i=1}^k\gamma_i\}, \ \ 1\leq k\leq n-1.\tag{4.1}\\
&\sum_{i=1}^n\alpha_i=r_1+\sum_{i=1}^{n-1}\beta_i.\tag{4.2}
\end{align}\end{thm}
{\bf Proof. } Use the Lemma \ref{lem:4.1} for the cases where all the invariant factors are powers of $p$.  $\Box$ 

\

Associating Lemma \ref{lem:4.2}, for counting sublattices in ${\mathcal S}_\alpha(p)$, we need to consider the admissible triples $(r_1, \Lambda^*_{n-1}, b_1^\star)$ satisfying the equations (4.1) and (4.2). 
	
For a partition $\beta=(\beta_1,\cdots,\beta_{n-1})$, fix an isomorphism of finite abelian groups $$\kappa:\Lambda_{n-1}/\Lambda_{n-1}^*\rightarrow{\mathbb Z}/p^{\beta_1}{\mathbb Z}\times\cdots\times{\mathbb Z}/p^{\beta_{n-1}}{\mathbb Z}.$$
For $b_1^\star\in\Lambda_{n-1}/\Lambda_{n-1}^*$, denote $\kappa(b_1^\star)=(c_1,\cdots,c_{n-1})$. Then, the coefficient matrix of $\Lambda^*_{n-1}(b_1^\star)$ with a suitable basis of $\Lambda_{n-1}$ is
$$M^\star=\begin{pmatrix}
	c_1\ c_2\cdots \ \  c_{n-1}   \\
\text{diag}\{p^{\beta_1},\cdots,p^{\beta_{n-1}}\}\\
\end{pmatrix}.$$ 

For $1\leq i\leq n-1$, define $$\delta_i=\delta_i(b_1^\star)=\begin{cases}\beta_i,&c_i\equiv 0 (\bmod\ \ p^{\beta_i})\\\text{ord}_p(c_i),&c_i\not\equiv 0 (\bmod\ \ p^{\beta_i})\end{cases}.$$
And, one has that $0\leq \delta_i\leq\beta_i, 1\leq i\leq n-1$.

For $1\leq k\leq n-1$, set $$\eta_k=\eta_k(\beta,\delta)=\min\{\beta_k-\beta_i+\delta_i,\delta_j\ |\ 1\leq i<k,k\leq j\leq n-1\},$$ and we take $\beta_0=\eta_0=0$ for convenience.
Then, $\eta_0\leq\eta_1\leq\cdots\leq\eta_{n-1}$. 

Recall that we assume $\Lambda_{n-1}^*(b^\star)\in\mathcal S_\gamma(p)$.
Considering the invariant factors of  $M^\star$, we have that $D_k(M^\star)=p^{\sum_{i=1}^{k-1}\beta_i+\eta_k}$ and
$$\sum_{i=1}^k\gamma_k=\sum_{i=1}^{k-1}\beta_i+\eta_k.$$
	
Denote 
$$k_0=\min\{1\leq k\leq n\ |\ r_1<\eta_k\},$$
i.e., $\eta_{k_0-1}\leq r_1<\eta_{k_0}$. $(4.1)$ and $(4.2)$ in Theorem \ref{thm:4.1} become
$$\sum_{i=1}^k\alpha_i=\sum_{i=1}^{k-1}\beta_i+\eta_k,\ \ 1\leq k< k_0$$ and
$$\sum_{i=1}^k\alpha_i=\sum_{i=1}^{k-1}\beta_i+r_1,\ \ k_0\leq k\leq n.$$
Making differences of the above sums, one has that
\begin{equation}\alpha_k=\begin{cases}\eta_1,& k=1\\\beta_{k-1}+\eta_k-\eta_{k-1},&2\leq k\leq k_0-1\\\beta_{k_0-1}+r_1-\eta_{k_0-1},&k=k_0\\\beta_{k-1}, &k_0+1\leq k\leq n\end{cases}.\tag{4.3}\end{equation}

On the other hand, we define $T(r_1,\alpha,\beta)$ is the set of all admissible $(n-1)$-tuples $\delta$(not need to be a partition) such that $0\leq\delta_i\leq\beta_i$ and$(4.3)$ holds. It is easy to see that for a fixed $\delta$, the number of admissible $b_1^\star\in\Lambda_{n-1}/\Lambda_{n-1}^*$ is equal to 
\begin{equation}\prod_{1\leq i\leq n-1, \delta_i<\beta_i}(p^{\beta_i-\delta_i}-p^{\beta_i-\delta_i-1}),\tag{4.4}\end{equation} which is an integral polynomial in $p$, denoted by $\tau(\beta,\delta)$.
	
By Lemma \ref{lem:4.2} and Theorem \ref{thm:4.1}, we obtain
\begin{thm}\label{thm:4.2} Using above notions, for a partition $\alpha\in{\mathcal P}(n,k)$,
\begin{equation}f_\alpha(p)=\sum_{\stackrel{0\leq r_1\leq k}{\beta\in{\mathcal P}(n-1,k-r_1)}}\left(f_\beta(p)\cdot\sum_{\delta\in T(r_1,\alpha,\beta)}\tau(\beta,\delta)\right).\ \ \Box\tag{4.5}\end{equation}
\end{thm}	
	
\
	
{\bf Remark. }Although the set $T(r_1,\alpha,\beta)$ is not clear enough to present $f_\alpha(p)$ explicitly, the formula may be useful to design an more effective iterative algorithm to compute $f_\alpha(p)$
than brute exhausting all the elements in $\text{HNF}_n(p^k)$.

At least, we can obtain the following property of $f_\alpha(p)$,
\begin{cor}\label{cor:4.1} Assume that $n\geq 1$. For any partition $\alpha\in{\mathcal P}(n,k)$, $f_\alpha(p)$ is an integral polynomial in $p$ with coefficients independent of $p$, i.e., there exists a polynomial $g_\alpha(T)\in\ZZ[T]$ such that for any prime $p$, $f_\alpha(p)=g_{\alpha}(p)$ holds.
\end{cor}	
{\bf Proof. }We induct on $n$.	When $n=1$, it is trivial. Assume that $f_\beta(p)$ are integral polynomials in $p$ with coefficients independent of $p$ for all $\beta\in{\mathcal P}(n-1,r)$.

For any $\delta\in T(r_1,\alpha,\beta)$, we have seen that $\tau(\beta,\delta)$ are integral polynomials in $p$ with coefficients independent of $p$ in $(4.4)$. The set $T(r_1,\alpha,\beta)$ is also independent of $p$. Hence, $f_\alpha(p)$ is an integral polynomial with coefficients in $p$ by the identity $(4.5)$ and the induction on $f_\beta(p)$. $\Box$

\

\section{Co-cyclic sublattices and their number}

In this section, we give an explicit formula on the number of co-cyclic sublattice with a fixed index $m$. First, we see that all co-cyclic sublattices with index $m$ belong to the unimodular equivalent class ${\mathcal S}(1,\cdots,1,m)$.

\begin{thm}\label{thm:5.1}All co-cyclic sublattices in an $n$-dimensional lattice $\Lambda$ with the same index are unimodular equivalent.\end{thm}
{\bf Proof. } Let $\Lambda^*$ is a sublattice of $\Lambda$. We see that if $\Lambda^*\in{\mathcal S}(d_1,\dots,d_n)$,
$$\Lambda/\Lambda^*\simeq(\ZZ/d_1\ZZ)\times\cdots\times(\ZZ/d_n\ZZ)$$
as finite abelian groups with $d_1\mid\cdots\mid d_n$.

Hence, $\Lambda^*$ is co-cyclic with index $m$ if and only if $d_1=\cdots=d_{n-1}=1$ and $d_n=m$, i.e., $\Lambda^*$ belongs to the equivalent class ${\mathcal S}(\underbrace{1,\cdots,1}_{n-1},m)$. $\Box$

The following Theorem gives the formula of the number $f(1,\cdots,1,p^r)$ of co-cyclic sublattice with the power-of-prime index $p^r$.

\begin{thm}\label{thm:5.2} For $n\geq 1$, $p^r$ is a power of prime $p$ with $r\geq0$,
	$$f(\underbrace{1,\cdots,1}_{n-1},p^r)=\begin{cases} 1  &r=0 \\p^{(n-1)(r-1)}\frac{p^n-1}{p-1}  &r\geq1\end{cases}.$$
\end{thm}	
{\bf Proof. }Let $\Lambda$ be an $n$-dimensional lattice with a basis $\{b_1,\cdots,b_n\}$. Observer that a sublattice $\Lambda^*$ of $\Lambda$ with coefficient matrix $M$ is co-cyclic if and only if $D_{n-1}(M)=1$.

The result is trivial when $n=1$. By induction on $n$, we assume the result is true until $n-1$.
	
Set $H=\begin{pmatrix}
		p^{r_1} & h_{12}\cdots h_{1n} \\
		0 & H' \\
	\end{pmatrix}$ is the coefficient matrix of $\Lambda^*$ and $H''=\begin{pmatrix}
		h_{12}\cdots h_{1n}  \\
		H' \\
	\end{pmatrix}
	$. 
By Lemma \ref{lem:4.1}, $\Lambda^*$ is co-cyclic if and only if the following identity holds
$$\gcd(p^{r_1}D_{n-2}(H'),D_{n-1}(H''))=1.$$  Let $b^\star=\sum_{i=2}^nh_{1i}b_i$ be the vector in the lattice $\Lambda_{n-1}$ spanned by $\{b_2,\cdots,b_n\}$ and $\Lambda^*_{n-1}$ the sublattice of $\Lambda_{n-1}$ with coefficient matrix $H'$. There are three cases can occur.
	
{\bf Case 1. } $r_1=0$. We must have that $D_{n-2}(H')=1$ or $D_{n-1}(H'')=1$. 

If  $D_{n-2}(H')=1$ , we have that  $\Lambda^*_{n-1}$ is a co-cyclic sublattice in $\Lambda_{n-1}$.  And, if $D_{n-1}(H'')=1$, we have that $\Lambda^*_{n-1}(b_1^\star)=\Lambda_{n-1}$ and $\Lambda^*_{n-1}$ is also co-cyclic in $\Lambda_{n-1}$. Hence, in this case, $\Lambda^*$ is co-cyclic in $\Lambda$ if and only if
$\Lambda^*_{n-1}$ is co-cyclic in $\Lambda_{n-1}$ with the same index $p^r$. And, any the choice of
$b_1^\star\in\Lambda_{n-1}/\Lambda^*_{n-1}$ is admissible. By Lemma \ref{lem:4.2}, the number of co-cyclic sublattices in this case is $p^r\cdot f(\underbrace{1,\cdots,1}_{n-2},p^r)$ in total.
	
{\bf Case 2. } $0<r_1<r$.  We have that $D_{n-1}(H'')=1$ and $\Lambda^*_{n-1}(b^\star)=\Lambda_{n-1}$. 

Then, $\Lambda^*_{n-1}$ is co-cyclic in $\Lambda_{n-1}$ with index $p^{r-r_1}$. The number of admissible $b_1^\star\in\Lambda_{n-1}/\Lambda^*_{n-1}$ is equal to the number of generators in the cyclic group $\Lambda_{n-1}/\Lambda^*_{n-1}$ with order $p^{r-r_1}$, which is $\varphi(p^{r-r_1})=p^{r-r_1-1}(p-1)$.
	
{\bf Case 3. } $r_1=r$. $\Lambda^*_{n-1}=\Lambda_{n-1}$ and $\Lambda^*$ is the unique sublattice corresponding to $\begin{pmatrix}
		p^r & 0\\
		0 & I_{n-1} \\
	\end{pmatrix}$.
	
Putting them together and using the induction on $n$, we have
\begin{align*}
	&f(\underbrace{1,\cdots,1}_{n-1},p^r)\\
	=&p^r\cdot f(\underbrace{1,\cdots,1}_{n-2},p^r)+\sum_{r_1=1}^{r-1}p^{r-r_1-1}(p-1)\cdot f(\underbrace{1,\cdots,1}_{n-2},p^{r-r_1})+1\\
	=&p^{(n-2)(r-1)+r}\frac{p^{n-1}-1}{p-1}+\sum_{r_1=1}^{r-1} p^{(n-1)(r-r_1-1)}(p^{n-1}-1)+1\\
		=&\sum_{i=0}^{n-2}p^{(n-1)r-i}+p^{(n-1)(r-1)}=p^{(n-1)r}+p^{(n-1)r-1}+\cdots+p^{(n-1)r-(n-1)}\\
		=&p^{(n-1)(r-1)}(p^{n-1}+\cdots+1)=p^{(n-1)(r-1)}\frac{p^n-1}{p-1}.  \ \ \Box
	\end{align*}

Notice that $f(1,\cdots,1,p^r)$ is an $n$-termed polynomial in $p$ with the leading term $p^{(n-1)r}$.
In fact, we apply $(4.5)$ in Theorem \ref{thm:4.2} for the simplest case indeed. 

\

When $m=p_1\cdots p_l$ is square-free, one has 
\begin{cor}\label{cor:5.1}For $n\geq 1$ and square-free $m$ with $m=p_1\cdots p_l$, we have
	$$f(\underbrace{1,\cdots,1}_{n-1},m)=\prod_{k=1}^l\frac{p_k^n-1}{p_k-1}.$$
\end{cor}
{\bf Proof. } By Theorem 5.2, 
$$f(1,\cdots,1,p)=\frac{p^n-1}{p-1}.$$
Using the multiplicativity in Corollary \ref{cor:2.2},
$$f(1,\cdots,1,m)=\prod_{k=1}^l\frac{p_k^n-1}{p_k-1}$$
when $m=p_1\cdots p_l$ is square-free. $\Box$

{\bf Remark.} We see that when $m$ is square-free, $f(1,\cdots,1,m)=f_n(m)$, i.e., all sublattices with index $m$ are unimodular equivalent and $g_n(m)=1$(Also see the Remark 4.4 in \cite{Z}).
 
\
 
For general case, one has
\begin{cor}\label{cor:5.2}For $n\geq 1$ and $m=p_1^{r_1}\cdots p_l^{r_l}$ where $p_1,\cdots, p_l$ are distinct primes, we have
$$f(\underbrace{1,\cdots,1}_{n-1},m)=\left(\frac{m}{r(m)}\right)^{n-1}\sigma_1(r(m)^{n-1}),$$
where $$r(m)=p_1\cdots p_l$$ is called the radical of $m$ and $\sigma_1(k)=\sum_{0<d\mid k}d$ is the sum of all positive divisors of an integer $k$. \end{cor}
{\bf Proof. }It is well-known that $\sigma_1$ is multiplicative. Also by the multiplicativity of $f$,
\begin{align*}&f(1,\cdots,1,m)=\prod_{k=1}^lf(1,\cdots,1,p_k^{r_k})\\
=&\prod_{k=1}^lp_k^{(n-1)(r_k-1)}(p_k^{n-1}+\cdots+1)=\frac{(p_1^{r_1}\cdots p_l^{r_l})^{n-1}}{{(p_1\cdots p_l)^{n-1}}}\prod_{k=1}^l\sigma_1(p_k^{n-1})\\
=&\frac{(p_1^{r_1}\cdots p_l^{r_l})^{n-1}}{{(p_1\cdots p_l)^{n-1}}}\sigma_1((p_1\cdots p_l)^{n-1})\\
=&\left(\frac{m}{r(m)}\right)^{n-1}\sigma_1(r(m)^{n-1}). \ \ \ \ \Box \end{align*}

\
	
In the last, we point out an interesting phenomenon that the co-cyclic sublattices take over the main part of all sublattices in an $n$-dimensional lattice.

When $n=2$, in Theorem \ref{thm:3.1}, we have seen that $f(1,p^r)=p^r+p^{r-1}$ and it is the two leading terms of $f_2(p^r)$ as polynomials in $p$. 

Assume that $n>2,r>0$. By (1.2), the number $f_n(p^r)$ of sublattices with index $p^r$ is $\prod_{k=1}^{n-1}\frac{p^{r+k}-1}{p^k-1}$. 
Although $(p^{r+k}-1)$ is not always divisible by $(p^k-1)$, a direct computation shows that
\begin{align*}&f_n(p^r)=\frac{\prod_{k=1}^{n-1}(p^{r+k}-1)}{\prod_{k=1}^{n-1}(p^k-1)}\\
	=&\frac{p^{\frac{n(n-1)}{2}+(n-1)r}-p^{\frac{n(n-1)}{2}+(n-2)r-1}+\cdots}{p^{\frac{n(n-1)}{2}}-p^{\frac{n(n-1)}{2}-1}+\cdots}\\
=&p^{(n-1)r}+p^{(n-1)r-1}+\cdots.\end{align*} 

Then, as polynomials in $p$, both $f_n(p^r)$ and $f(1,\cdots,1,p^r)$ have the same two leading terms $p^{(n-1)r}+p^{(n-1)r-1}$ at least. Particularly,
$$\frac{f(1,\cdots,1,p^r)}{f_n(p^r)}\rightarrow 1 \text{ as } p\rightarrow\infty.$$ 

Hence, almost sublattices of an $n$-dimensional lattice
with index $p^r$ are co-cyclic for large enough $p$. In other words,  for an integer $m$, the larger prime factors does $m$ have, the more possible is that a sublattice with index $m$ is co-cyclic. This gives another aspect of the fact that most sublattices are co-cyclic mentioned in \cite{NS}. 
	
{\bf Remark.} In \cite{NS}, they claim that “most” integer lattices are co-cyclic for sufficiently large dimension. In our aspect, also for a fixed dimension, co-cyclic sublattices take over the main part of sublattices with index $m$, especially when all the prime factors of $m$ are large enough.

\end{document}